\documentclass[12pt]{article}
\usepackage{amsmath,amsfonts}
\begin{document}

\baselineskip 20pt
\date{}
\title{An affirmative answer to a problem of Cater\footnote{2010 Mathematics Subject Classification: Primary 26A24; 26A48 Key words: Increasing function, derivative.}}

\author{Arthur A.~Danielyan}

\maketitle

\begin{abstract}
 
\noindent  Does there exist an increasing absolutely continuous function, $f: [0,1] \rightarrow \mathbb R$ such that $\{x: f'(x)=0\}$ is both countable and dense? This problem was proposed by F.S. Cater about two decades ago. We give an affirmative answer to the problem.

\end{abstract}

\begin{section}{Introduction and the main result.}

The following problem was proposed by F. S. Cater \cite{Cater}.

  \vspace{0.25 cm}

 {\bf Problem 1.} {\it Does there exist an increasing absolutely continuous function, $f: [0,1] \rightarrow \mathbb R$ such that $\{x: f'(x)=0\}$ is both countable and dense?} 

  \vspace{0.25 cm}
  
The question is referring to an increasing function, but $f$ must be strictly increasing because if $f$ is not strictly increasing (but just increasing) then the set $\{x: f'(x)=0\}$ contains an interval and cannot be countable. Thus, obviously, in Problem 1 the word ``increasing" has the meaning of ``strictly increasing".
 
We give an affirmative answer to Problem 1, by proving the following theorem. 

 \vspace{0.25 cm} 
 
 {\bf Theorem 1.} {\it There exists an 
 increasing absolutely continuous function, $f: [0,1] \rightarrow \mathbb R$ such that $\{x: f'(x)=0\}$ is both countable and dense.} 

\vspace{0.25 cm}

Recall that the lower derivative ${\underline f'}(x)$ of a function $f$ at $x_0$ is defined by
$$ {\underline f'}(x_0) =  \liminf_{x \rightarrow x_0}\frac{f(x)-f(x_0)}{x-x_0} $$
The upper derivative is defined similarly (see e.g. \cite{Nat}, p. 158).  

 The following result shows that 
the conclusion of Theorem 1 is precise and that it cannot be extended for the case of lower derivative. 

 \vspace{0.25 cm} 
 
 {\bf Theorem 2.} {\it There does not exist an 
 increasing continuous function, $f: [0,1] \rightarrow \mathbb R$ such that $\{x: {\underline f'}(x)=0\}$ is both countable and dense.} 

\vspace{0.25 cm}

Note that Theorem 2 concerns any increasing continuous functions (and not just increasing absolutely continuous functions). 

The proof of Theorem 2 will be given in another publication.

\end{section}

\begin{section}{An auxiliary result.}

We will use the following known theorem due to V. M. Tzodiks \cite{Tzo}.

    \vspace{0.25 cm}

{\bf Theorem A.}\ {\it Let $E$ be an $F_{\sigma \delta}$ set of measure zero and let $N$ be a $G_{\delta}$ set, 
such that both $E$, $N$ are on $\mathbb R$, and $N \supset E$. Then there exists an increasing continuous function $F(x)$
such that: (1) $F'(x)= + \infty$ on $E$;  (2)  ${\underline F'}(x) < + \infty$ for $x \not \in E$; and (3)   $F'(x)$ exists and is finite
for $x \not \in N$.}
 
\vspace{0.25 cm}

\end{section}

\begin{section}{Proofs.}
{\it Proof of Theorem 1.} We construct an increasing absolutely continuous function $f$ on $[0,1]$ such that $\{x: f'(x)=0\}$ is both countable and dense. 

Let $R$ be a countable and dense subset of numbers in the open interval $(0,1)$. For example, one can take as $R$  the set of all rational numbers in  $(0,1)$.

 Since $R$ is an $F_\sigma$ set, it is also 
an $F_{\sigma  \delta}$ set. We can easily construct a $G_\delta$ set $G$ of measure zero such that 
$R \subset G \subset [0,1]$. 

We apply Theorem A taking $E=R$ and $N=G$. This gives an increasing continuous
 function $F$ such that $F'(x)= + \infty$ on $R$;  
${\underline F'}(x) < + \infty$ for $x \not \in R$; and $F'(x)$ exists and is finite
for $x \not \in G$. 
(Theorem A provides the function $F$ defined on entire $\mathbb R$, but we only consider its restriction on the interval $[0,1]$.)

In particular, $F'(x)$ exists and is finite a.e. on $[0,1]$. Since $F$ is increasing, $F'(x) \ge 0$ whenever  $F'(x)$  
 exists. 
 
 Without losing the generality, we may assume that $F(x)$ does not decrease the distance of any two points, as one can simply replace  $F(x)$ by  $x+F(x)$, if needed. 
 Indeed, obviously the function $x+F(x)$ too possesses the properties of $F(x)$ listed above.

Since $F$ is increasing and continuous on $[0,1]$, $F$ maps $[0,1]$ onto the interval $[F(0),F(1)]$.  
Without losing the generality we may assume that $[F(0),F(1)]=[0,1]$. 
Indeed, one can just replace
$F(x)$ by $$F_1(x)=\frac{1}{F(1)-F(0)}[F(x) - F(0)], $$ 
and the latter function inherits all other properties of the former.

Since $F(x)$ does not decrease the distance of any two points, $F(x)-F(0)$ does the same, and thus for $0 \leq x_1< x_2 \leq 1$, 

$$F_1(x_2) - F_1(x_1) \geq \frac{1}{F(1)-F(0)}(x_2-x_1).$$

The last inequality directly implies that the inverse $F_1^{-1}$ of $F_1$ is a Lipschitz 1 function with constant $[F(1)-F(0)].$
Thus, $F_1^{-1}$ is absolutely continuous.

Since $F_1$ is increasing and continuous, as well as $R$ is countable and dense on $[0,1]$,
the image set $F_1(R)$ of $R$ is countable and dense on $[0,1]$. 

Let $f=F_1^{-1}.$ Then $f$ is increasing and absolutely continuous on $[0,1]$.

Because $F'_1(x)=+ \infty$ on $R$ and $f$ is the inverse of $F_1$,  $f'(x)=0$ on the set $F(R)$.
Recall that ${\underline {F'_1}(x)} < + \infty$ for $x \not \in R$; this implies that $f'(x)$ is not zero for  $x \not \in F(R)$.
Thus the zero set of $f'(x)$ is $F(R)$. 

The proof of Theorem 1 is over. 

\vspace{0.25 cm} 

{\bf Acknowledgement.} The author wishes to thank Vilmos Totik for a helpful 
remark which simplified the original proof of Theorem 1.

\end{section}

\begin{minipage}[t]{6.5cm}
Arthur A. Danielyan\\
Department of Mathematics\\ 
and Statistics\\
University of South Florida\\
Tampa, Florida 33620\\
USA\\
{\small e-mail: adaniely@usf.edu}\\
\end{minipage}


\begin{thebibliography}{99}

\bibitem{Cater} F. S. Cater, A countable and dense zero set, Real Analysis Exchange,   {\bf 33(2)}, 2007/2008, pp. 483 - 484. 

\bibitem{Nat} I. P. Natanson, Theory of functions of a real variable, V. 2, Ungar, 1960.

\bibitem{Tzo} V. M. Tsodiks, On the sets of points where the derivative is finite or infinite correspondingly, Dokl. Akad. Nauk SSSR
(N.S.) {\bf 114}, 1957, 1174-1176.




\end{thebibliography}
\end{document}